\documentclass[11pt]{amsart}
\usepackage{latexsym, amsfonts, amsmath, amssymb}

\newlength{\standardunitlength}
\newcommand{\id}{\mathrm{id}}
\newcommand{\D}{\Delta}
\newcommand{\ot}{\otimes}

\newcommand{\lan}{\langle}
\newcommand{\ran}{\rangle}

\newcommand{\eps}{\varepsilon}

\def\CM{{\mathcal M}}
\def\CC{{\mathcal C}}

\def\BC{{\mathbb C}}

\def\Hom{\mathrm{Hom}}

\newtheorem{thm}{Theorem}[section]
\newtheorem{df}[thm]{Definition}
\newtheorem{remark}[thm]{Remark}

\newtheorem{cor}[thm]{Corollary}
\newtheorem{prop}[thm]{Proposition}
\newtheorem{lem}[thm]{Lemma}

\newcommand{\num}{\refstepcounter{thm}}

\newtheorem{lemma}[thm]{Lemma}

\newtheorem{theorem}[thm]{Theorem}

\newcommand{\RR}{\mathbb{R}}
\newcommand{\QQ}{\mathbb{Q}}

\newcommand{\ZZ}{\mathbb{Z}}

\newcommand{\C}{\mathbb{C}}

\newcommand{\Ng}{N_g}

\newcommand{\alt}{\mathrm{Alt}}
\newcommand{\sym}{\mathrm{Sym}}
\newcommand{\Alt}{\mathrm{Alt}}
\newcommand{\Sym}{\mathrm{Sym}}

\begin{document}

\begin{center} \title [Schur Indicators] {\bf Frobenius-Schur
Indicators for Subgroups and the Drinfel'd Double of Weyl Groups}
\end{center}

\author{Robert Guralnick}
\address{University of Southern California \\
Los Angeles, CA 90089-2532}
\email{guralnic@usc.edu}

\author{Susan Montgomery}
\address{University of Southern California \\
Los Angeles, CA 90089-2532}
\email{smontgom@usc.edu}



\date{\today}

\thanks{The authors were supported by NSF grants DMS 0140578 and
   DMS 0401399.}

\begin{abstract} If $G$ is any finite group and $k$ is a field, there
is a natural construction of a Hopf algebra over $k$ associated to
$G$, the Drinfel'd double $D(G)$. We prove that if $G$ is any
finite real reflection group, with Drinfel'd double $D(G)$ over an
algebraically closed field $k$ of characteristic not $2$, then
every simple $D(G)$-module has Frobenius-Schur indicator +1. This
generalizes the classical results for modules over the group
itself.  We also prove some new results about Weyl groups. In
particular, we prove that any abelian subgroup is inverted by some
involution. Also, if $E$ is any elementary abelian $2$-subgroup of
the Weyl group $W$, then all representations of $C_W(E)$ are
defined over $\QQ$.
\end{abstract}

\maketitle


\section{Introduction} \label{intro}

  There has been considerable interest lately in using Frobenius-Schur
indicators
in studying semisimple Hopf algebras, as they are one of only a few useful
invariants
for a Hopf algebra $H$. For example indicators (and their higher analogs)
have been used
in classifying Hopf algebras themselves \cite{K}
\cite{NS1}; in studying possible dimensions of the
representations of $H$ \cite{KSZ1}; and in studying
the prime divisors of the exponent of $H$
\cite{KSZ2}. Moreover, the indicator is invariant under
equivalence of monoidal categories
\cite{MaN}. Another motivation comes from conformal field theory;
see work of Bantay \cite{B1}
\cite{B2}. Thus it is important to compute the values of the
indicator for more examples.

Let $k$ be an algebraically closed field of characteristic not
$2$, $H$ an  involutive Hopf algebra (i.e. the antipode
has order $2$)  and $V$ a finite dimensional irreducible
$H$-module.   We show  (generalizing the case of
group algebras and semisimple Hopf algebras) that one can define
the Schur indicator  of $V$ and it has the following properties:
 $\nu(V)\ne 0$ if and only if $V$ is self
dual, $\nu(V)=+1$ if and  only if $V$ admits a non-degenerate
$H$-invariant symmetric bilinear form and $\nu(V)=-1$ if and only
$V$ admits a non-degenerate $H$-invariant skew-symmetric bilinear
form. If $G$ is a finite group, then self duality is equivalent to
$\chi(x)=\chi(x^{-1})$ for all $x \in G$ where $\chi$ is the
character of $V$.  In particular, if $k = \BC$, this is saying
that the character is real valued (and all representations are
self dual if and only if $x$ and $x^{-1}$ are conjugate in $G$ for
all $x$).  Also, $\nu(V)=+1$ if and only if $V$ is defined over
$\RR$.

It was shown in \cite{LM} that the classical theorem of Frobenius and Schur,
giving a formula to
compute the indicator of a simple module for a finite group $G$ extends to
semisimple Hopf algebras. This
fact was shown a bit earlier for the
special case of Kac algebras over $\BC$ in \cite{FGSV}.
In  this paper, we show that the indicator is well defined
whenever the antipode
of $H$ satisfies $S^2 = \id$.  Moreover,  the indicator
has the above characterization in terms of bilinear
forms.  The formula for the indicator is no longer valid (this is clear already
for modular group algebras).

The case when $\nu(V) = +1$ for all simple $G$-modules $V$
has been of particular interest; such
groups are called {\it totally orthogonal} in \cite{GW}.
This terminology seems
suitable for Hopf
algebras as well, in view of the existence of the bilinear
forms described above.
In particular it was
known classically that if $G$ is any finite real reflection
group, then all $\nu(V) =
+1$ \cite{H}; other
examples are given in \cite{GW}.

For Hopf algebras, there has also been some interest in when a Hopf algebra is
``totally orthogonal'',
in particular when $H = D(G)$, the Drinfel'd double of a finite group.
It is proved in
\cite{KMM} that if
$G$ is a dihedral group or the symmetric group $\Sym_n$ and $k = \BC$,
then $D(G)$ is totally
orthogonal.
More generally, \cite{KMM} studies the indicator for Hopf algebras
which are abelian extensions,
since in this case the irreducible representations are known.

A natural question to ask is what happens for the other Weyl groups, since
their group algebras in any characteristic not $2$ has all
indicators for irreducible
modules equal to $+1$.
The  main result we prove for Drinfel'd doubles is:

   \begin{theorem} \label{main}  Let $G$ be a finite real
reflection group and $D(G)$
   its Drinfel'd double
over the algebraically closed field $k$ of characteristic not $2$.
Then $D(G)$ is totally orthogonal; that is, every simple $D(G)$-module
has Frobenius-Schur indicator $+1$.
   \end{theorem}

  Recall that a real reflection group is a subgroup of $GL(V)$
with $V$ finite dimensional
  over $\RR$  that is generated by
  reflections.  Every such is a direct product of irreducible real
reflection groups.  The finite
  irreducible real reflection groups are well known (cf \cite{H}):
aside from
  some examples in dimension at most $4$, they are the Weyl groups
and for dimension
  greater than $8$, they are just the classical Weyl groups.

We also prove some new general results for real
reflection groups and Weyl groups
which give more than just our main Theorem and are of independent interest.
We prove these results by considering first the classical Weyl groups
(that is, those of type A, B and D).   For the smaller
reflection groups, we use
computer calculations -- some of which were quite simple
and some more involved.
We refer the reader to  \cite{breuer} and \cite{gg} for more details.

We should point out that it is false in general that if $G$ is totally
orthogonal, then the same is true for $D(G)$: we give an
example of a group of order 32 with the property that
all of its representations
are defined over $\QQ$ but such that $D(G)$ has an irreducible
module with indicator 0 (see Section \ref{example}).
Thus a proof of Theorem \ref{main}
will have to involve the specific properties of the reflection groups.
Indeed,
we shall see that the classical Weyl groups satisfy some stronger
properties than
all real reflection groups (or even all Weyl groups).

   We now discuss some of our results about Weyl groups and real
reflection groups,
   If $g \in G$, let $C_g$ denote the centralizer of $g$ and $N_g$
   the normalizer of $\{g,g^{-1}\}$.

   \begin{theorem} \label{inversion}  Let $G$ be a finite real
reflection group and let $A$ be an abelian subgroup of $G$.
\begin{enumerate}
\item There exists an
involution $t \in G$
   that acts as inversion on $A$.
\item If $G$ is a Weyl group and $r$ is relatively prime to $|A|$,
then there exists $x \in G$ with $xax^{-1}=a^r$ for all $a \in A$.
\end{enumerate}
   \end{theorem}

   The first part of the theorem
 is a generalization of the well known property for
Weyl groups that every
   element is a product of two involutions (cf. \cite[\S8]{springera}).
This also implies
   that all representations of $C_g$ are self dual, but we show more:

   \begin{theorem} \label{rational}  Let $G$ be a real reflection group.
    Let $E$ be an elementary abelian $2$-subgroup of $G$ and let $k$
   be an algebraically closed field of characteristic $p \ne 2$.
   \begin{enumerate}
   \item Every element of $C_G(E)$ is a product of $2$ involutions and
    $kC_G(E)$ is totally orthogonal; and
   \item If $G$ is a  Weyl group,
   then every irreducible complex representation of
   $C_G(E)$ is defined over $\mathbb{Q}$.
   \end{enumerate}
   \end{theorem}

   Let $L$ be a field.
   Recall that the field of definition for a module $V$ for $LG$ is the
   smallest field $K$ such that $V = U \otimes_K L $ with $U$ a $KG$-module
   (of course, $K$ must contain all character values and in positive
characteristic,
   this is sufficient as well for $V$ semisimple).

   We prove the second assertion in Theorem \ref{rational}
in this paper for the
   classical Weyl groups.  It is proved for the exceptional Weyl groups
   in \cite{gg}.   In particular, this generalizes the well known result that
   all representations of Weyl groups are defined over $\QQ$.  Our proof is
   independent of this (we do use the relatively elementary
result for symmetric groups).

   We also prove:

   \begin{theorem} \label{Ng}  Let $G$ be a real reflection group
   and $k$ an algebraically closed field of characteristic $p \ne 2$.
   Let $g \in G$.
   \begin{enumerate}
   \item If $G$ is a classical Weyl group or is a real reflection group of
   dimension at most $3$, then all Frobenius-Schur indicators for
   irreducible $kN_g$-modules are $+1$.
   \item For  all other real reflection groups, all Frobenius-Schur
   indicators for irreducible $kN_g$-modules are non-negative.
   \end{enumerate}
   \end{theorem}

   One computes easily that for the Weyl groups of type $E$ and $F$ and for
   the real reflection group $H_4$, there exists an element of order $3$
   such that $N_g$ is not totally orthogonal.  Thus,
   one cannot expect a uniform proof.

   We will show that the two  previous theorems are enough to show that
   $D(G)$ is totally orthogonal for $G$ a finite real reflection group.  A curious
   consequence of the result on $D(G)$ gives a formula for the finite reflection groups
   that seems new:

   \begin{cor} \label{formula}
   Let $G$ be a finite reflection group.  Let $t$ be the number
   of elements of $G$ whose square is trivial.  Then
   $$
   t^2 = \sum_{g \in G} \sum_{\chi \in \mathrm{Irr}(C_g)} \chi(1).
   $$
   \end{cor}

   Indeed, this formula holds for a finite group if and only if $D(G)$ is totally
   orthogonal.   This follows since $t^2$ is the trace of the antipode of $D(G)$ and when all
   indicators are $+1$, this is equal to the right hand side (see \cite{LM}, \cite{KMM}).

   The paper is organized as follows.   In Section \ref{NScond} we review and
   extend the results from
   \cite{KMM} which we need, and give a refinement of the criterion for positivity
   of Schur indicators in \cite{KMM}.  We also generalize the result to 
   any characteristic other than $2$
   and do not assume that $D(G)$ is semisimple. In order to do this, we first show
   that the indicator is well defined when $H$ is not necessarily semi-simple as long
   as $S^2= \id$. In section \ref{Schur Indices}, we review some basic properties
   of Schur indices of finite groups.

   In Section \ref{dpsection} we consider the case of direct products
   and we also show that it is sufficient to prove Theorem \ref{main}
   in  characteristic $0$. In Section \ref{abelian subgroups}, we consider
   normalizers of abelian subgroups of Weyl groups.  In particular, 
    we prove Theorem \ref{inversion} and the results about centralizers of involutions.
    In Section \ref{realgroupsec}, we show that certain groups are split over $\RR$.
    In Section
   \ref{typeB} we consider the classical Weyl groups.
   In Section
   \ref{general} we consider the remaining real reflection groups.
   For these small cases, we obtain the results by computation (some quite easy and
   some more involved).  We also prove Corollary \ref{formula}.
   Finally in
   Section \ref{example} we discuss the example showing total orthogonality
   does not always
   lift from $G$ to $D(G)$.


\section{Indicators for Involutive Hopf Algebras and Drinfel'd Doubles}\label{NScond}

We recall a few basic facts about modules over a Hopf algebra H over a field $k$.
We assume throughout that $S^2 = \id$, for the antipode $S$ of
$H$. This is always the case by \cite{LR} for semisimple
Hopf algebras in characteristic zero and so is also the case for Hopf algebras
in positive characteristic which may be defined by reducing semisimple Hopf algebras
modulo $p$ (eg, group algebras, Drinfe'ld doubles of finite groups and
and bismash products of groups) . Masuoka has pointed
out to us that it is
also true for bicrossed products of groups, using descent arguments
as in \cite{Mat}.

Let $\CC = {_H{\CM}}$ denote the category of
finite-dimensional left $H$-modules. Let $V_0 = kv_0\in \CC$
denote the {\it trivial} module, that is $h\cdot v_0 = \eps(h)
v_0$. If $V,W \in \CC$, then $V\ot W \in \CC$, via $h\cdot (v\ot
w) := \D (h)\cdot (v\ot w) =
  \sum h_1\cdot v\ot h_2\cdot w$, for all $h\in H, v\in V, w\in W$, and
  also $\Hom_k(V,W) \in \CC$
  via $(h\cdot f)(v) := \sum h_1\cdot (f(Sh_2\cdot v)) \in W$.
As a special case, the dual module $V^* \in \CC$ also, via
$(h\cdot f)(v) = f(Sh\cdot v)$ for all $h\in H, f\in H^*$ and
$v\in V$. Since $S^2=id$, $V^{**} \cong V$.

For any $V\in \CC$, let $V^H:= \{v \in V \;|\;h\cdot v = \eps(h) v
\;\forall h\in H\}$, the {\it $H$- invariants} of $V$.

In \cite{LM}, the indicator was defined for semisimple Hopf
algebras in terms of a formula involving integrals. This
definition does not work here, since in the non-semisimple case,
there does not exist an integral $\Lambda \in \int _H$ with
$\epsilon (\Lambda)=1$. Instead we use the conclusion of
\cite[Theorem 3.1, (1) and (2)]{LM} as our definition. This agrees
with the classical definition for finite groups in characteristic
$p$.

\begin{df} {\rm Let $H$ be a Hopf algebra over an algebraically
closed field of characteristic not 2 such that $S^2 = \id$, and
let $V$ be a finite-dimensional irreducible (left) $H$-module. We
define a function $\nu: V\to \{0, 1, -1\}$, called the {\it
Frobenius-Schur indicator} of $V$, as follows:

(1) If $V$ is not self-dual, set $\nu(V) = 0$.

(2) If $V^* \cong V$ and $V$ admits a non-degenerate $H$-invariant
symmetric bilinear form,  set $\nu(V) = +1$.

(3) If $V^* \cong V$ and  $V$ admits a non-degenerate
$H$-invariant skew-symmetric form, set $\nu(V) = -1$}.
\end{df}

We prove here that in fact, any irreducible $H$-module $V$ is of
precisely of type (1), (2), or (3). We first need an elementary
lemma, most of which is well-known.

\begin{lem}\label{rep} Let $V,W\in \CC$. Then as left $H$-modules,
\begin{enumerate}
\item $\Hom(V,W^*) \cong \Hom(V\ot W, k)$ and $\Hom_H (V,W)=
\Hom(V,W)^H$.
\item  Assume $k$ is algebraically closed and $V$ is irreducible. If
$V^* \cong V$, then $\Hom_H(V\ot V, k)$ is one-dimensional, and if
$V^* \ncong V$, then $\Hom_H(V\ot V, k)=0.$
\item Assume $k$ is algebraically closed, $V$ is irreducible, $V^*
\cong V$, and choose $0\neq \langle-,- \rangle \in \Hom_H(V\ot V,
k)= ((V\ot V)^*)^H$. Then
$$ \sum  \langle h_1 \cdot v, h_2 \cdot w \rangle=\eps(h) \langle v,w \rangle.$$
\item  The form $\langle-,- \rangle \in \Hom_H(V\ot V, k)$ $\iff$ the
antipode $S$ is the adjoint with respect to the form.
\end{enumerate}
\end{lem}

\begin{proof} (1) The first part   is straightforward, using the standard map \hfil\break
$\Hom(V,W^*) \to \Hom(V\ot W, k)$ given by $f \mapsto \tilde f$,
where $\tilde f(v\ot w):= f(v)(w)$. The second fact follows from
the definitions.

(2) follows from (1) using $W = V$ and Schur's Lemma.

(3) The displayed equation is simply saying that $\langle-,-
\rangle$ is an $H$-module map. The fact that the form may be
considered as an invariant of $(V\ot V)^*$ follows from (1).

(4) Finally, we consider the relationship of the form $\langle-,-
\rangle$ to $S$. We use the property that $S$ is a ``generalized
inverse'' map on $H$, that is for any $h\in H$ with $\D(h) = \sum
h_1\ot h_2$, $S(h_1)h_2 = \eps(h)$; we also use the counit
property of $\eps$ that $\sum \eps(h_1)h_2 = h$, for all $h\in H$.

 First assume that $\langle-,- \rangle$ is $H$-invariant. Then for all $h\in H,v,w\in V$,
\begin{eqnarray*}
\langle h\cdot v,w \rangle &=& \sum \langle h_1\cdot v,\eps(h_2)w
\rangle =\sum \langle h_1\cdot v,h_2 S(h_3) \cdot w \rangle \\
&=& \sum  \langle (h_1)_1 \cdot v,(h_1)_2 \cdot(S(h_2)\cdot w)
\rangle =\sum \eps(h_1) \langle v,S(h_2)\cdot w \rangle \\
&=& \langle
v,Sh \cdot w \rangle.
\end{eqnarray*}
Conversely, assume that $\langle h\cdot v,w \rangle =\langle v,Sh
\cdot w \rangle$ for all $h,v,w$. Then $$\sum \lan h_1\cdot v, h_2
\cdot w \ran = \sum \lan v, S(h_1)h_2 \cdot w\ran = \lan
v,\eps(h)w \ran =\eps(h) \lan v,w \ran.$$ and so the form is
$H$-invariant.
\end{proof}

From (2), every $H$-invariant form on $V\cong V^*$ is a scalar
multiple of the given one $0\neq \langle-,- \rangle$. Given a form
$\langle-,- \rangle \in V\cong V^*$, we may also consider its
``twist'' $\tau$, that is the form $[v,w] := \tau( \langle v,w
\rangle) = \langle w, v \rangle.$

The argument in the next lemma is taken from \cite[p 18]{KSZ2}.

\begin{lem} \label{twist} If $\langle-,- \rangle$ is $H$-invariant,
then so is its twist $[v,w]  = \langle w, v \rangle.$
Consequently $\langle w, v \rangle = \alpha \langle w, v \rangle$,
for all $v, w$, where $\alpha= \pm 1$ is independent of the choice
of  $\langle-,-\rangle$ and depends only on $V$.
\end{lem}

\begin{proof} Using $S^2 = \id$ and Lemma \ref{rep},
$$[h \cdot v,w] = [S^2(h) \cdot v, w] = \langle w, S^2(h)\cdot v \rangle=
\langle S(h)\cdot w, v \rangle = [v, S(h)\cdot w], $$ for any $v,
w \in V$, $h\in H.$

Since $((V\ot V)^*)^H$ is one dimensional, $[-,-] = \alpha
\langle-,-\rangle$ for some $\alpha \in k$. Repeating the argument
gives $\alpha^2 = 1$, and so $\alpha = \pm 1$. It follows that
$\alpha$ is independent of the choice of  $\langle-,-\rangle$ and
depends only on $V$.
\end{proof}

The theorem is now clear from our arguments above.

\begin{thm} \label{Sind} Let $H$ be a Hopf algebra such that $S^2 = id$,
over the algebraically closed field $k$ of characteristic $p \ne 2$
and let $V$ be a finite dimensional  irreducible (left)
$H$-module. Then $V$ falls into case (1), (2), or (3) of the
definition, and so has a Frobenius-Schur indicator. Moreover,
\begin{enumerate}
\item  $V^* \cong V$ $\iff$ $\nu(V)\neq 0$ $\iff$ $(V\ot V)^H$ is
one-dimensional; and
\item  if $\nu(V) \neq 0$ and $f$ is a
non-zero invariant in $V\ot V$, then $\tau(f) = \nu(V)f$, where
$\tau$ is the twist map above.
\end{enumerate}
\end{thm}

We now consider Drinfel'd doubles.
Throughout the remainder of the section, let $k$ be an algebraically closed field of
characteristic $p$ with $p \ne 2$.

For a finite group $G$, let $kG$ be the
group algebra and $(kG)^*$ its dual Hopf algebra,
with dual bases $\{g \in G\}$ and $\{p_g|g\in G\}$. The Drinfel'd double $D(G)$ of $G$ is
simply $(kG)^*\ot kG$ as a vector space. As an algebra it is the smash product
$(kG)^* \# kG$, where
$kG$ acts on $(kG)^*$ via $g\cdot p_h = p_{ghg^{-1}}$, for all $g,h \in G$.
The coalgebra structure is given by the tensor product of coalgebras, and
the antipode $S$ is given by
$$ S(p_x \# a)=p_{ax^{-1}a^{-1}} \# a^{-1}.$$
One sees directly that $S^2 = \id$.

The description of the simple
$D(G)$-modules is well-known \cite{DPR}, \cite[Section 2]{Ma},
and in fact depends only on the algebra structure of $D(G)$. More precisely, for any
$g\in G$ let $C_g$
denote the centralizer of $g$ in $G$, and choose one $g$ in each conjugacy class of $G$.
For each irreducible $C_g$-module $V$, let $\hat V =  V_{C_g}^G$ be the induced module;
then $\hat V$ is a $kG$-module, and it becomes a $D(G)$-module via
\num
\begin{equation}
(p_h \#y) \cdot [x \ot  v] = \delta_{x g x^{-1},h} [yx \ot  v].
\end{equation}
for $h,x \in G$ and $v \in V$.

Then $\hat V$ is an irreducible $D(G)$-module, and all irreducible $D(G)$-modules
arise in this way.
Thus the simple $D(G)$-modules $\hat V$ are indexed by pairs $(g,V)$ where $g \in G$
ranges
over a set of conjugacy class representatives and $V$ is an irreducible $C_g$-module.

We want to  refine the result \cite[Theorem 5.5]{KMM} about
Schur indicators of the simple $H$-modules.  We also extend the results
to the case of positive characteristic without the assumption that $D(G)$
is semisimple.  Thus, our proof will both be more general and
independent of \cite{KMM}.

As in \cite{KMM}, we will show that to compute $\nu(\hat V)$,
one does not need to
induce $V$ up to $G$, but only up to $N_g$, the normalizer of the set
$\{g,g^{-1}\}$ in $G$,
and then use the usual indicator for the representations of the groups
$C_g$ and $N_g$:

\begin{theorem}  \label{localind} Let $k$ be an algebraically
closed field of characteristic not $2$, $G$  a finite group, $D(G)$
the Drinfel'd double of $G$ over $k$ and  $\hat{V}$ an irreducible
$D(G)$-module corresponding to the irreducible $C_g$-module $V$.
\begin{enumerate}
\item If $g^2=1$, then $\nu(\hat{V}) = \nu(V)$.
\item If $g^2 \ne 1$, then $\nu(\hat{V}) = \nu(V_{C_g}^{\Ng}) - \nu(V)$.
\end{enumerate}
\end{theorem}

We will prove this below in a series of results which give a bit more
information.

For $V$ a $C_g$-module, we will denote the conjugate module to $V$ by $V^x$; that is,
$V^x = V$ as a $k$-space, but for any $g\in C_g$ and $v\in V^x$,
$g\cdot v = xgx^{-1}\cdot v$.

If $h \in G$, let  $I(h)$ be the $k$-span of
$\{ p_x \#y | y \in G, x \in h^G\}$.  Then $I(h)$ is a $2$-sided ideal
of $H$ and $H$ is the direct sum of the $I(h)$ as $h$ ranges over a set
of conjugacy class representatives.
  Moreover, $I(1) \cong kG$. Since $\epsilon(x)=0$ for $x \in I(h), h \ne 1$,
we see that  the $D(G)$ invariants on a module $W$ are precisely
the $G$-fixed points on $I(1)W$.

We first observe:

\begin{lemma} \label{step1}  Let $g \in G$ and $V$ be an irreducible
$kC_g$-module. If $g$ and $g^{-1}$ are not conjugate, then $\nu(\hat{V})=0$.
\end{lemma}

\begin{proof}   We note that $I(h) \hat{V}=0$ if and only if $h$
is not conjugate to  $g$ while $I(h) {\hat{V}}^* = 0$ if and only
if $h$ is not conjugate to $g^{-1}$. So if $g$ and $g^{-1}$ are not conjugate,
then clearly $\hat{V}$ is not self dual.

\end{proof}

\begin{theorem} \label{localind2}  Let $g \in G$ and $V$ be an irreducible
$kC_g$-module. Set $C=C_g$ and $N=N_g$. Let $x \in G$ with $xgx^{-1}=g^{-1}$.
\begin{enumerate}
\item $\nu(\hat{V}) \ne 0$ if and only if $V \cong V^x$ as $C_g$-modules.
\item If $g^2=1$, then $\nu(\hat{V})=\nu(V)$.
\item If $g^2 \ne 1$, then $\nu(\hat{V})=\nu(V_C^N)-\nu(V)$.
\end{enumerate}
\end{theorem}

\begin{proof}  We may assume that $x=1$ if $g^2=1$. Note that $x \in N$.

Let $S$ be a set of left coset representatives for $N$ in $G$.
Then $T:=S \cup Sx$ is a set of left coset representatives for
$C$ in $G$. So $\hat{V} =
\sum_{t \in T} t \otimes V$.
Set $V_t:=t \otimes V$.

Set $X = \hat{V} \otimes \hat{V}$ and let $\tau$ be the natural
involution on $X$.

It follows easily from the definitions that
$Y:= I(1)X = \sum_{t \in T} V_t \otimes V_{tx}$.
Thus the $D(G)$-invariants on $X$ are precisely the
$G$-fixed points on $Y$.  Note that $Y \cong (V \otimes V_x)_C^N$.
By Frobenius reciprocity, $G$ has fixed points on $Y$ if and
only if $C$ has fixed points on $V \otimes V_x$.  Since $V_x \cong V^x$
as $kC$-modules, we see that $D(G)$ has invariants on $X$ if and only
if $V \cong V^*$ as $C$-modules.  This proves (1).

Note that Frobenius reciprocity gives a bit more. It tells us that
if $f$ is a nonzero $C$-fixed point on $V \otimes V_x$, then
$F:=\sum_{t \in T} tf$ is a nonzero fixed point of $G$ on $Y$
(and so a nonzero $D(G)$-invariant).

Suppose that $g^2 =1$.  Then $V \otimes V$ is  $\tau$-invariant and
$\tau$ commutes with the action of $C$ on $V \otimes V$.
So if $f$ is a nonzero $C$-fixed element, then $\tau(f)= \nu(V)f$
and $\tau(F)=\nu(V) F$, whence $\nu(V)=\nu(\hat{V})$.  This proves
(2).

Now assume that $g^2 \ne 1$. If $V^x$ and $V^*$ are non-isomorphic
$C$-modules, then it is straightforward to compute 
that $\nu(V)=\nu(V_C^N)$ and so (3) holds in this case.

So finally, we assume also that $V^x$ and $V^*$ are isomorphic
as $C$-modules.

Let $U = (V \otimes V_x) \oplus (V_x \otimes V)$.  Note that $U$ is
an $N$-module and is $\tau$-invariant.  Also, it is clear
that $Y \cong U_N^G$.   Thus, arguing as above, we deduce
that there is a nonzero $N$-fixed point $f$ in $U$
and  $\tau(f) = \nu(\hat{V})f$. We have now reduced the
computation of $\nu(\hat{V})$ to a computation in $N$ and we
proceed to prove (3).

Let $W=V_C^N$ and note that $U$ embeds in $W \otimes W$.
Also, the action of $\tau$ on $U$ extends to that of
$W \otimes W$ in the obvious manner (with $\tau$ the natural
involution on $W \otimes W$).   The result now follows by
the group theoretic result, Lemma \ref{localindgroup} below.
\end{proof}

\begin{lemma} \label{localindgroup}  Let $k$ be an algebraically
closed field of characteristic $p \ne 2$.  Let $C$ be a subgroup
of index $2$ in the finite group $N$.  Fix an element $x \in N \setminus{C}$
and let $V$ be an irreducible $kC$-module with $V^* \cong V^x$.  Set $W=V_C^N
=V \oplus xV$.  Let $X = W \otimes W$ and let $\tau$ be the linear map
sending $w_1 \otimes w_2$ to $w_2 \otimes w_1$.   Let $Y = (V \otimes xV) \oplus
(xV \otimes V)$, a $kN$-submodule of $X$.   The fixed points of $N$ on $Y$
are $1$-dimensional and are $\tau$-invariant.  Let $f$ be a generator
for this space.   Then $\tau(f) = (\nu(W) - \nu(V))f$, where
$\nu(W)$ is the Frobenius-Schur indicator of $W$ as a $kN$-module and
$\nu(V)$ is the Frobenius-Schur indicator of $V$ as a $kC$-module.
\end{lemma}

\begin{proof}  Note that $X = Y \oplus Z$ where $Z=V \otimes V \oplus xV \otimes xV$.
Clearly $Y$ and $Z$ are both invariant under $N$ and $\tau$.  Note that
$C$ has a $1$-dimensional fixed space on $V \otimes xV \cong V \otimes V^*$.
Let $w$ be a generator for this space.  Then $w, xw$ is a basis for
the fixed space of $C$ on $Y$ and so $f:=w + xw$ is a basis for the fixed
space of $N$ on $Y$.   Since $\tau$ commutes with $N$ on $W$, $\tau f = \gamma f$
where $\gamma = \pm 1$.  So we only need show that
$\gamma =\nu(W) - \nu(V)$.

 We consider various cases.\\

\noindent Case 1.   $V$ is not self dual as a $kC$-module.\\

Then $W$ is an irreducible self dual $kN$-module, whence
the $N$-invariants on $W \otimes W$ are $1$-dimensional and so generated by
$f$.  Thus,  $\gamma = \nu(W)$ and since $\nu(V)=0$, the result holds.\\

\noindent Case 2.   $V$ is self dual.   By Frobenius reciprocity,
$W = U_1 \oplus U_2$ where the $U_i$ are nonisomorphic $kN$-modules
such that $U_i \cong V$ as $kC$-modules.   Note that the space of
$C$-invariants on $Y$ and $Z$ are each $2$-dimensional and the space
of $N$-invariants on $Y$ and $Z$ are each $1$-dimensional.  The
invariants of $C$ on $Z$ are contained in the $\nu(V)$-eigenspace
of $\tau$.   Since $\tau$ interchanges $V \otimes xV$ and $xV \otimes V$,
we see that $\tau$ has two eigenvalues on the invariants of
$C$ on $Y$.  Thus, $\tau$ has a $3$-dimensional $\nu(V)$-eigenspace
and $1$-dimensional $-\nu(V)$-eigenspace on the $C$-invariants on $X$.

First suppose that each $U_i$ is self dual as a $kN$-module.
Then $\nu(U_i) = \nu(V)$ for each $i$ and so $\nu(W)=2 \nu(V)$ and
we need to show that $\gamma = \nu(V)$ in this case.  We can also
view $X = \oplus U_i \otimes U_j$ and we see that the $N$-invariants
are contained in $U_1 \otimes U_1 \oplus U_2 \otimes U_2$ and so
we see that $\tau$ acts on this space as $\nu(V) = \nu(W) - \nu(V)$.

Finally, suppose that each $U_i$ is not self dual.  Thus, $\nu(W)=0$.
We need to show that $\gamma = - \nu(V)$ in this case.
On $Y$, $\tau$ acts via $\nu(V)$ on the $C$-fixed points and so on
the $N$-fixed points.   Since
$W$ is self dual (as $V$ is), this implies that $U_2=U_1^*$.
Thus, we see that $\tau$ has both eigenvalues on the $2$-dimensional
space of $N$-fixed points on $W$, whence $\tau(f) = - \nu(V) f$.
\end{proof}

It is straightforward to see that this implies the following theorem,
which is a  refinement
of \cite[Corollary 6.2]{KMM} if $k = \C$ but is new if $k$
has positive characteristic $p$ dividing $|G|$.

\begin{prop} \label{real} Let $G$ be a finite group.
Then $\nu(\hat{V}) = 1$ for all irreducible $D(G)$-modules $\hat V$ if and only if:
\begin{enumerate}
\item When $g^2=1$,
 $\nu(V) = 1$ for all irreducible $C_g$-modules $V$;
\item When $g^2 \ne 1$, $|\Ng:C_g|=2$ and for any irreducible $C_g$-module $V$,
\begin{enumerate}
   \item $V^x \cong V^*$ for all $x \in \Ng \setminus{C_g}$, where $V^x$ is the
   conjugate module to $V$ ;
   \item if $V$ is not self dual, then $V_{C_g}^{\Ng}$ has Schur
          indicator $+1$,
   \item if $\nu(V)=+1$,
      then each constituent of $V_{C_g}^{\Ng}$ has Schur indicator $+1$.
   \item  if $\nu(V)=-1$,
          then neither constituent of $V_{C_g}^{\Ng}$ is self dual.

\end{enumerate}
   \end{enumerate}
   \end{prop}

   \begin{proof}  If $g^2=1$, this is clear.  This also implies that
   every element of $\langle g \rangle$ is conjugate to its inverse.

    So consider the case where $g^2 \ne 1$ and fix $x \in \Ng$ with
   $xgx^{-1} = g^{-1}$. Recall that $V^x = V$ as a $k$-space, but with new action
   $g\cdot_x v = xgx^{-1}\cdot v$, for all $g\in C_g$, $v\in V^x$.

   We have already seen that $V^x \cong V^*$ if and only if all indicators
   for $D(G)$ are nonzero.  Note that the condition that $V^x \cong V^*$
   implies that $V_{C_g}^{N_g}$ is self dual.  So we assume that this is
   the case.

   Suppose that $V$ is not self dual.  Then
   $\nu(\hat{V}) = \nu(V_{C_g}^{\Ng})$.   So we may assume that
   $V$ is self dual and $V \cong V^x$.

   Now it follows by Frobenius reciprocity
   that the induced module is a direct sum of two irreducible
   modules. Note that the sum of the characters of the
   two irreducible modules vanishes outside $C_g$ (because the sum
   of the modules is an induced module).  In particular, one is self dual if
   and only if the other one is.  So if neither is self dual, then
   $\nu(\hat{V}) = -\nu(V)$ and the theorem holds.

   If they are both self dual, then the Schur indicators for each of
   the two modules must be equal to $\nu(V)$ and so $\nu(\hat{V}) = \nu(V)$.
   This completes the proof.

\end{proof}

In particular, this shows:

\begin{cor} \label{to for Ng} Let $G$ be a finite group and
$k$ an algebraically closed field of characteristic not $2$.
Assume that $g$ is conjugate to $g^{-1}$ for all $g \in G$.
If $x \in N_g \setminus{C_g}$ implies that $V \cong V^x$
for all irreducible $kC_g$-modules and $\nu(W)=1$ for
all irreducible $kN_g$-modules, then $D(G)$ is totally orthogonal.
\end{cor}

The previous results also imply:

   \begin{cor} \label{selfdual}
    Let $G$ be a finite group. Then $\nu(\hat{V}) = \pm 1$ for all irreducible
$D(G)$-modules $\hat V$ if and only if:
   \begin{enumerate}
    \item when $g^2=1$, every irreducible module of $C_g$ is self dual;
    \item when $g^2 \ne 1$, $|\Ng:C_g|=2$ and for any irreducible $C_g$-module $V$, we
have $V^x \cong V^*$ for all (any) $x \in \Ng \setminus{C_g}$.
    \end{enumerate}
    \end{cor}

\begin{remark} {\rm G. Mason has also observed that Corollary \ref{selfdual} for
$k = \BC$ can
be obtained using the methods of \cite{KMM} (private communication).}
\end{remark}

  \begin{remark} {\rm The condition that $V^x \cong V^*$ for all $\BC[C_g]$ modules
   $V$ with $x \in \Ng \setminus{C_g}$ is equivalent to saying that
   $x$ inverts all conjugacy classes of $C_g$ (i.e. if $y \in C_g$, then
   $y^x$ is conjugate to $y^{-1}$ in $C_g$).}
\end{remark}


\section{Schur Indices I}  \label{Schur Indices}

We recall some basic facts about Schur indices.  See \cite{Feit}
for details.  Let $G$ be a finite group and $k$ a field of characteristic
$0$.  Let $V$ be an absolutely irreducible $G$-module over the
algebraic closure of $k$ with character $\chi$.  Let $k'$ be the extension field of
$k$ generated by the values of $\chi$.
The Schur index $m_{\chi}$ is the smallest positive integer such
that there is a $k'G$-module $W$ with character $m_{\chi} \chi$.

The basic fact that we will require about this is (cf. \cite[Lemma 9.1]{Feit}):

\begin{lemma}  \label{basic schur index}
Suppose $W$ is a $kG$-module with character $\theta$.
If $\chi$ is an irreducible character of $G$ with values in $k$, then
$m_\chi$ divides $(\theta,\chi)$.
\end{lemma}

This has the following easy consequence:

\begin{lemma} \label{index 2} Let $N$ be a subgroup of index $2$ in the finite
group $G$.   Let $k$ be a subfield of $\BC$.
Let $V$ be an irreducible $\BC N$-module with character $\chi$ having
values in $k$.
Set $W=V_N^G$.
\begin{enumerate}
\item If each irreducible constituent of $W$ is defined over
$k$, then   $V$ is defined over $k$.
\item If $V$ is defined over $k$ and each irreducible constituent of
$W$ has character values in $k$, then they are defined over $k$.
\end{enumerate}
\end{lemma}

\begin{proof}
The key point is that $W$ is multiplicity free.  Indeed, by Frobenius
reciprocity, either $W$ is irreducible and $W_N$ is a direct sum of two
non-isomorphic modules ($V$ and $V^x$ for $x \in G \setminus{N}$) or
$W$ is the direct sum of two non-isomorphic $G$-modules which are
isomorphic as $N$-modules.

(1) If $W$ is not irreducible and is defined over $K$, then each constituent
is isomorphic to $V$ as an $N$-module, whence $V$ is defined over $K$.

If $W$ is irreducible, then by Frobenius reciprocity,
$1=(\chi_N^G,\chi_N^G)_G=(\chi,\chi_N^G)$.  Thus, the previous result
implies that $m_\chi=1$ as required.

(2)  If $W$ is irreducible, then since it is induced from a module
defined over $k$, it is defined over $k$.  If $W$ is not irreducible,
then it is the sum of two non-isomorphic irreducible modules each with
character in $k$.  By the previous result, the Schur index of
each constituent is $1$, whence they are each defined over $k$.
\end{proof}

An immediate consequence is the following:

\begin{cor} \label{index 2 cor}
Let $N$ be a subgroup of index $2$ in the finite
group $G$.   Let $k$ be a field of characteristic $0$.
Suppose that all characters of $N$ and $G$ are defined over $k$.
Then $k$ is a splitting field for $N$ if and only if it is a splitting
field for $G$.
\end{cor}

Let $C$ be a subgroup of index $2$ in the finite group $N$.
Pick $x \in N \setminus{C}$ and assume that $x$ inverts
all conjugacy classes of $C$.  We consider two conditions:

\begin{enumerate}
\item[DC1]  If $V$ is an irreducible $C$-module with $\nu(V) = 0$,
then $V_C^N$ is irreducible with $\nu(V_C^N)=+1$.  Moreover,
if $\nu(V) = \pm 1$, then $V$ has two (distinct) extensions
to $N$, both with Schur indicator $+1$ if $\nu(V)=+1$ and $0$
if $\nu(V)=-1$ (this can be summarized by saying that
$\nu(V_C^N) - \nu(V_C) = + 1$ for all irreducible $C$-modules).
\item[DC2]$\nu(W)=+1$ for all irreducible $N$-modules.
\end{enumerate}

It is an easy exercise to show that DC2 implies that DC1.
We show that these properties behave nicely with respect to
direct product.

\begin{lemma} \label{dp1}
Let $N_i, 1 \le i \le t$ be finite groups with subgroups $C_i$
of index $2$. Fix $x_i \in N_i \setminus{C_i}$.
Assume that $x_i$ inverts all conjugacy classes of $C_i$.
Let $N = N_1 \times \ldots \times N_t$ and $C= C_1 \times \ldots \times C_t$.
Set $H = \langle C, x \rangle$ where $x=x_1 \cdots x_t$.
\begin{enumerate}
\item $x$ inverts all conjugacy classes in $C$.
\item If $C_i < N_i$ satisfies $\mathrm{DC1}$ or
$\mathrm{DC2}$ for each $i$, then so does $C < H$.
\end{enumerate}
\end{lemma}

\begin{proof}  The first statement is clear.  This shows that
if $V$ is an irreducible $C$-module and is not self dual, then
$V^x \cong V^*$ is not isomorphic to $V$.

We now prove (2). By induction,  we may reduce to the case
that $t=2$.   If $\nu(W)=+1$ for all irreducible $N$-modules
(i.e. DC2 holds), then
by Corollary \ref{index 2 cor}, it suffices to show that all characters of
$H$ are real valued.  Clearly, every element of $C$ is conjugate to
its inverse in $H$.  Suppose that $y:=(y_1,y_2) \in H$ and $y_i$ is not
in $C_i$.  Then $y_i$ is conjugate to its inverse via an element
of $C_i$ since all representations of $N_i$ are self dual.
Then $y^g=y^{-1}$ for some $g \in C$ whence $y$ is conjugate to
its inverse in $H$.

We now show that DC1 extends as well.

Suppose that $V$ is an irreducible $C$-module.
Write $V=V_1 \otimes V_2$.  So $\nu(V)=\nu(V_1)\nu(V_2)$.

If $\nu(V_1)=\nu(V_2)=1$, then both $V_i$ extend to irreducible $N_i$-modules
all with Schur indicator $+1$ and so therefore the same is true for
each constituent of $V_C^H$.

If $\nu(V_1)=\nu(V_2)=-1$, then
$\nu(V)=+1$ and so $V$ is defined over $\RR$.   By hypothesis,
$U_i:={V_i}_{C_i}^{N_i}$ is a direct sum of two irreducible modules which
are not self dual for $N_i$.   Let $X_i$ and $Y_i$ be the two irreducible constituents.
Since $U_i$ is self dual, $Y_i=X_i^*$ and since the character of $U_i$ vanishes
outside $C_i$, it follows that the character of $X_i$ (and of $Y_i)$ is purely
imaginary on $N_i \setminus{C_i}$.   Thus, the character of $X_i$ on $H$ is real
valued.  By Lemma \ref{index 2}, this implies that both $X_i$ and $Y_i$ are
defined over $\RR$.

If $\nu(V_1) = - \nu(V_2) \ne 0$, then $\nu(V)=-1$ and we need to show that
an extension of $V$ to $H$ is not self dual.  Assume that $\nu(V_1)=-1$.
Then we can choose $c_i \in C_i$ such that the character of $x_ic_i$ is nonzero
and purely imaginary for $i=1$ and real for $i=2$ (the character cannot vanish
outside $C_i$ by Frobenius reciprocity).   Then the character value of $(c_1x_1,c_2x_2)$
on any extension of $V$ is purely imaginary and nonzero (being the product of the
character values for the extensions of $V_1$ and $V_2$) as required.

So we may assume that $\nu(V)=\nu(V_1)=0$.   Then $W:=V_C^H$ is irreducible and self dual.
By hypothesis,  $W_H^N = V_C^N$.       So $U_1$ is irreducible
as an $N_1$-module and $\nu(U_i)=1$.   Note that $C_1$ has a $1$-dimensional fixed
space on both $\mathrm{Sym}^2(U_1)$ and $\wedge^2(U_1)$ and that $H$ is trivial on
$\mathrm{Sym}^2(U_1)$ and $x_1$ acts as $-1$ on $\wedge^2(U_1)$.

Note that $W_H^N=V_C^N=U_1 \otimes U_2$.  If $\nu(V_2) \ge 0$, then each constituent
of $U_2$ has Schur indicator $+1$ and so $W_H^N$ is defined over $\RR$, whence
by Lemma \ref{index 2}, $W$ is.

Suppose that $\nu(V_2) = -1$.   Then $C_2$ has a fixed point on $\wedge^2(V_2)$.
Now $V_2$ can be viewed as a non-self dual module for $N_2$ (in two distinct ways)
and $W = U_1 \otimes V_2$.  Thus,
$x_2$ acts on $\wedge^2(V_2)$ by $-1$.  So we see that $x$ has a fixed point
on $\wedge^2(U_1) \otimes \wedge^2(V_2) \le \mathrm{Sym}^2(U_1 \otimes V_2)$.
Thus, $\nu(W)=+1$ as required.
\end{proof}


   \section{Some Reductions} \label{dpsection}

   Let $G$ be a finite group, and consider the following three properties of $G$:

   \begin{enumerate}
   \item[(R1)]  Each irreducible $D(G)$-module has Schur indicator $\pm 1$;
   \item[(R2)]  Each irreducible $D(G)$-module has Schur indicator $ +1$;
   \item[(R3)]  Each irreducible $D(G)$-module has Schur indicator $ +1$
      and for each $g \in G$, every $\Ng$-module has Schur indicator $+1$.
   \end{enumerate}

   It is clear that (R3) implies (R2) implies (R1).  None of the
implications
   can be reversed.  For example, $G=Q_8$ satisfies (R1) but not (R2).
   The Weyl groups of types $F_4, E_6, E_7, E_8$ and the real reflection
    group $H_4$  satisfy (R2) but not (R3); this can be seen using Magma.  We will see
that the classical Weyl groups do satisfy (R3).

We describe the groups $G$ for which $D(G)$
satisfy  (R3)  using Proposition \ref{real}.

\begin{lemma} \label{R3}
Let $G$ be a finite group.
Then $D(G)$ (over $\BC$) satisfies $\mathrm{(R3)}$ if and only if for each
$g \in G$, $\nu(V)=+1$ for all irreducible $N_g$-modules $V$ and if
$x \in G$ inverts $g$, then $x$ inverts all $C_g$ conjugacy
classes.
\end{lemma}

   We first note:

   \begin{lemma} \label{semidirect}
   Suppose that  $G = N \rtimes H$, a semidirect product of a
   normal subgroup
   $N$ and a complement $H$. If $G$ satisfies {\rm (Ri)}
   with $1 \le i \le 3$, then so does $H$.
   \end{lemma}

   \begin{proof}  Let $h \in H$.  Then $C_H(h)$ is a homomorphic image of
   $C_h$ and similarly for $N_h$.  Thus, the modules considered in
   computing
   the Schur indicators have already arisen in the computation for $G$.
   \end{proof}

   \begin{prop}  \label{direct products} Let $G$ be a finite group.
   Assume that $G$ is a direct product of groups $G_j$.
   Then $G$ satisfies (Ri) if and only if each $G_j$ does.
   \end{prop}

   \begin{proof}  It suffices to consider a direct product of 2 groups.
   The previous lemma implies the forward implication.  So assume that
   each direct factor satisfies (Ri).

   So assume that $G= A \times B$ and $g=(a,b)$.  Then
   $C_g=C_a \times C_b$ and we can take $t=(u,v)$ inverting $g$.
   Then $t$ clearly inverts every class of $C_g$.  This shows that
   (R1) holds, using Theorem \ref{selfdual}.  So we are in the case of (R2) or (R3).
   Moreover, every $D(G)$-module has nonzero Schur indicator (since the result
   holds for (R1)).

   If $g^2=1$, then since $C_g$ is a direct product, the result is clear.

   So we may assume that $a^2 \ne 1$.  Let $V$ be an irreducible $C_g$
   module.  So $V =V_A \otimes V_B$ where $V_A$ and $V_B$ are irreducible
   $C_a$ and $C_b$-modules, respectively.  If $b^2=1$, then $\Ng = N_a
   \times C_b
   =N_a \times N_b$ and all irreducible modules are tensor products.
    Since
   $B$ satisfies (R2), all Schur indices of $C_b$-modules are $+1$.
    Since
   $A$ satisfies (R2) or (R3), it is trivial to see that all
   Schur indices for the $H$-modules occurring are $+1$ and that moreover
   if (R3) holds for $A$ and $B$, then every $\Ng$-module is real.

   So we may assume that $a^2 \ne 1 \ne b^2$.  Let $N=N_a \times N_b$.
   Now the result follows by Lemma \ref{dp1}.

   \end{proof}

We end this section by showing that being totally orthogonal in characteristic
zero implies the same result for odd   characteristic.
Let $k$ be an algebraically closed field of odd characteristic.
A result of Thompson \cite{thompson}  asserts that if
$V$ is an irreducible $kG$-module of nonzero Schur indicator, then
$V$ is a composition factor (of odd multiplicity) in the reduction mod $p$ of
an irreducible $\BC G$-module with the same Schur indicator as $V$.

\begin{theorem} \label{modular reduction}  Let $G$ be a finite group.
Let $k$ be an algebraically closed field of odd characteristic.
\begin{enumerate}
\item If all irreducible $\BC G$-modules have Schur indicator $+1$, the same
is true for all irreducible $kG$-modules.
\item If all irreducible $\BC G$-modules have Schur indicator at least $0$,
the same
is true for all irreducible $kG$-modules.
\item If all irreducible $\BC G$-modules are self dual, the same is true
for all irreducible $kG$-modules.
\end{enumerate}
\end{theorem}

\begin{proof} (3) is clear since the hypothesis implies all elements of
$g$ are conjugate to their inverse.

Now consider  (2).  By the result of Thompson mentioned above,
there are no irreducible $kG$-modules $V$ with $\nu(V)=-1$.
So (2) follows immediately.

Now (1) follows by (2) and (3).
\end{proof}

We use the previous result to show that $D(G)$ being totally orthogonal
in characteristic zero implies the same result for odd positive characteristic.

\begin{cor} \label{mod red for D(G)}
Let $G$ be a finite group and $k$ an algebraically
closed field of odd characteristic $p$.   Let $D(G)$
denote the Drinfel'd double
of $G$ over $\BC$ and $D_k(G)$ the Drinfel'd double of $G$ over $k$.
\begin{enumerate}
\item If $D(G)$ is totally orthogonal, then so is $D_k(G)$.
\item If $D(G)$ satisfies $\mathrm{(R3)}$, then does does $D_k(G)$.
\item If $D(G)$ satisfies $\mathrm{(R1)}$, then does does $D_k(G)$.
\end{enumerate}
\end{cor}

\begin{proof} Assume that $D(G)$ is totally orthogonal.  Note
that this implies that $G$ is totally orthogonal.
In particular, every element is conjugate to its inverse.  Fix $g \in G$.

If $g^2=1$,
then $C_g$ is totally orthogonal in characteristic $0$
and so by the previous result also over $k$.

Suppose that $g^2 \ne 1$.  Then, as noted above, $[N_g:C_g]=2$ and
$x \in N_g \setminus{C_g}$ implies that $x$ inverts
all conjugacy classes of $C_g$ (by the assumption for
characteristic $0$) and so $V^* \cong V^x$ for
all irreducible $kC_g$-modules $V$.

Let $V$ be an irreducible $kC_g$-module.

If $V$ is not self
dual, then $W:=V_{C_g}^{N_g}$ is self dual for $N_g$
and $\nu(\hat{V})=\nu(W)$. We claim that $\nu(W)=1$.

By Thompson's result, $W$ is contained in the reduction
of an irreducible $\BC N_g$-module $U$ with $\nu(U)=\nu(W)$
and the multiplicity of $W$ in the reduction of $U$ is odd.
Suppose that $U$ is irreducible as a $C_g$-module.
Since
$D(G)$ is totally orthogonal, it follows that
$1 =  \nu(U_{C_g}^{N_g}) -  \nu(U_{C_g}) = \nu(U_{C_g})=\nu(U_{N_g}=1$.
Since $\nu(U_{N_g})= \nu(W)$, the claim holds in this case.

If $U$ is not irreducible as a $C_g$-module, then $U=Y_{C_g}^{N_g}$ for some
irreducible $C_g$-module $Y$.  Then, $V$ is a constituent of the reduction
of $Y$ of odd multiplicity.  So if $Y$ is self dual, it follow that
$V$ is, a contradiction.  Thus, $V$ is not self dual and so $D(G)$
totally orthogonal implies that $\nu(W)=1$.

Suppose that $V$ is self dual and
so  $\nu(V) \ne 0$. Let $W$ be an irreducible $\BC C_g$-module
with $\nu(W)=\nu(V)$ such that $V$ occurs with odd multiplicity
in the reduction of $W$.   Then $W$ can be viewed as $N$-module.
Since $D(G)$ is totally orthogonal, either $\nu(W)=1$ both as
an $N_g$ and $C_g$ module or $\nu(W)=-1$ as a $C_g$-module and
$\nu(W)=0$ as an $N_g$-module.

First suppose that $\nu(W)=+1$ (and so also for $N_g$).
The reduction of $W$ has one of the extensions of $V$ occurring
with odd multiplicity (since this is true of $V$ in the reduction
of $W$ as a $C_g$-module).  Since $W$ is self dual, the composition
factors in the reduction must be closed under duals, whence
the $N_g$-module $V$ must be self dual and therefore must
satisfy $\nu(V)=+1$.

Suppose that $\nu(V)=-1$.   Since $\nu(U) \ne -1$ for any irreducible
$\BC N_g$-module, by Thompson's result, $\nu(V) \ne -1$ as a $kN_g$-module,
whence $\nu(V)=0$.

We have now shown that the necessary (and sufficient) conditions for $D(G)$
to be totally orthogonal over $k$ are satisfied.

If in addition, $D(G)$ satisfies (R3), we only need to verify that
$kN_g$ is totally orthogonal, but this is Theorem \ref{modular reduction}.

The proof for (R1) is analogous and we leave this as an exercise
to the reader (this is not used in the sequel).
\end{proof}


\section{Abelian Subgroups of Real Reflection Groups} \label{abelian subgroups}

We now prove Theorem \ref{inversion} for the classical Weyl
groups.  This also proves that every irreducible representation of
$C_g$ is self dual, but we will prove more later. In fact, we
prove a stronger result.  We first consider the case of symmetric
and alternating groups.

\begin{lemma} \label{abelian in sym}  Let $W=\Sym_n$ and let $J$
be an abelian subgroup of $W$.   Let $r$ be an integer prime to
$|J|$.  Then there exists $x \in W$ such that
$x^{-1}yx=y^r$ for all $y \in J$.   Moreover,
\begin{enumerate}
\item If  $r=-1$, then we can choose
$x$ to be an involution.
\item  If all nontrivial orbits of $J$ have even size, then we can choose $x \in \Alt_n$.
\item If $r=-1$ and all nontrivial orbits of $J$ have even size, we can choose $x$ to
be an involution in $\Alt_n$.
\end{enumerate}
\end{lemma}

\begin{proof}   We prove the first statement.  By induction, it suffices
to assume that $J$ is transitive and so acts regularly and has order $n$.
Let $X = \mathrm{Aut(J)}$.  Note that $H:=J.X$ embeds in $\Sym_n$ by letting $H$
act on the cosets of $X$.   Then $J$ maps onto a transitive abelian
subgroup of $\Sym_n$ and the result follows.

If $r=-1$, the same proof shows that  $x$ can be chosen to be an involution.

Suppose that all nontrivial orbits of $J$ have even order.  By induction,
we may reduce to the case that $J$ is transitive of even order.   If $J$
is not contained in $\Alt_n$, then we may replace
$x$ by $xs$ for some $s \in J$
to allow us to choose $x$ in either coset.  So $J \le \Alt_n$
which implies that
the Sylow $2$-subgroup of $J$ is not cyclic.

So we can write $J=C \times D \le S_c \times S_d$ with $C$ and $D$ of even
order and $n=cd$. Now choose $x_C$ and $x_D$
acting on $C$ and $D$ respectively in the appropriate
manner and take $x=x_Cx_D$.
We may assume that $c$ is even and so by induction,
$x_C \in A_c$ and so $x_C$ is also in $A_{n}$.
Since $c$ is even, $x_D$ is also in $A_{n}$
(whether or not $x_D \in A_d$).   Thus, $x \in \Alt_n$ as required.
Moreover, if $r=-1$, by construction $x^2=1$.
\end{proof}

We can now handle the classical Weyl groups.

\begin{theorem} \label{strong abelian}  Let $W$ be a
Weyl group, $J$ an abelian subgroup and $r$ an integer prime to
the order of $J$.  Then there exists $x \in W$ such that
$x^{-1}yx=y^r$ for all $y \in J$.  If $r=-1$, then we can choose
$x$ to be an involution.
\end{theorem}

\begin{proof} We first note that we may and do assume that $r$ is odd.
For if $J$ has even order, then certainly $r$ is odd.  If $J$ has odd
order, we may replace $r$ by $r + |J|$.
The case $W=\Sym_n$ follows by the previous result.

Next we consider $W$ of type $B_n$.  Let $W_0 \cong D_n$
be a subgroup of index $2$ in $W$.  We show that we can
take $x \in W_0$ with the appropriate properties.

We can view $W$ as the centralizer of a fixed point free
involution $z \in \Sym_{2n}$.  Note that $W_0=W \cap \Alt_{2n}$.

Let $J$ be an abelian subgroup of $W$.
There is no harm in assuming that
$z \in J$ (and so all orbits of $J$ have even order).
The result now follows by the previous result.

By the result for the symmetric
group, we can find $x \in \Alt_{2n}$ with $x^{-1}yx=y^r$ for
all $y \in J$.  In particular, $x$ commutes with $z$ and so
$x \in W \cap \Alt_{2n}=W_0$.   Moreover, if $r=-1$, we can
take $x$ to be an involution in $W_0$.  This completes the 
proof for the classical Weyl groups.  

If $W$ is exceptional, then Breuer \cite{breuer} has used GAP
to check that the result holds. 
\end{proof}

\begin{cor} Let $E$ be an elementary abelian $2$-subgroup of the
Weyl group $W$.  Then all characters of $C_W(E)$ are rational
valued.
\end{cor}

\begin{proof}  Let $g \in C_W(E)$ and let $r$ be an odd integer
prime to the order of $g$.  Then $J:=\langle E, g \rangle$ is
abelian and so there exists $x \in W$ with $xax^{-1}=a^r$ for all
$a \in J$.  Thus, $x \in C_W(E)$ and so $g$ is conjugate to any
power of $g$ prime to its order, whence all characters of
$C_W(E)$ are rational valued on $g$.
\end{proof}

For general real reflection groups, we observe:

\begin{theorem} \label{abelian} Let $W$ be a finite real
reflection group
and let $J$ be an abelian subgroup of $W$.
There exists an involution in $W$ that inverts $J$.
\end{theorem}

\begin{proof}   Clearly, we may assume that $W$ is irreducible.
If $W$ is a Weyl group, the result follows
by the previous theorem. If $W$ is dihedral or $W \cong H_3 \cong
A_5 \times \ZZ/2$, the result is clear.

The remaining groups are
$H_4$, $F_4$, $E_6$, $E_7$ and $E_8$.   In these cases, we use
computer calculations to verify the result (see  \cite{breuer} for details).
\end{proof}

\begin{cor}  Let $W$ be a real reflection group.
\begin{enumerate}
\item If $g^2=1$, then every element of $C_g$ is a product of two
involutions and in particular, all irreducible representations of
$C_g$ are self dual.
\item If $g^2 \ne 1$ and $x \in N_g \setminus{C_g}$, then $x$ inverts
all conjugacy classes of $C_g$.
\end{enumerate}
\end{cor}

\begin{proof}  Suppose that $g^2=1$ and $y \in C_W(g)$.
By the previous result, there exists an involution $x \in C_W(g)$
inverting $y$.  Then $y = x(xy)$ is a product of two involutions
in $C_W(g)$.  This implies that all irreducible representations
of $C_W(g)$ are self dual and proves (1).

Fix $x \in N_g$.
Let $y \in C_g$. By the previous result, there exists $x' \in N_g$
inverting $y$.  Thus, $x'=xc$ with $c \in C_g$.
Thus, $y^{-1} = y^{xc}$ and so the $C_g$ conjugacy class of
$y^x$ is the $C_g$ conjugacy class of $y^{-1}$ as asserted.
\end{proof}

Using Theorem \ref{selfdual},
this gives a fairly quick proof that the indicators
for the simple modules of
the Drin'feld double of a real reflection group
are at least nonzero.

We have seen that all characters of $C_W(E)$ are rational valued
for $E$ an elementary abelian $2$-group of the Weyl group $W$.
In fact, it is true all representations of $C_W(E)$ are defined
over $\QQ$.  We give the proof for the classical Weyl groups
below.  See \cite{gg} for the case of the exceptional Weyl groups.

\begin{theorem} \label{rat-inv}
Let $W$ be either a  Weyl group or an alternating group.
Let $E$ be an elementary abelian $2$-subgroup and if $W=\mathrm{Alt}_n$
assume that $E$ has no fixed points.
Then all representations of $C_W(E)$ are defined
over $\QQ$.
\end{theorem}

\begin{proof}  First suppose that $W=\mathrm{Sym}_n$.
If $E=1$, the result is clear.  If $1 \ne x \in E$ has
fixed points, then $C_W(E) \le C_W(x)$ is contained in
a direct product of two smaller symmetric groups and so
the centralizer will be the direct product of centralizers
in these smaller symmetric groups and the
result follows by induction.  So we may assume that all orbits
of $E$ are regular.  So $n=mt$ where $m=|E|$ and $t$ is the number
of orbits of $E$.  Then $C_W(E) \cong E \wr S_t$.  It is straightforward
to see that all representations of this group are defined over $\QQ$
(see the proof of Lemma \ref{real groups} for a more general result).

Now suppose that $E$ has no fixed points and $W=\mathrm{Alt}_n$.
Let $G=\mathrm{Sym}_n$.  If $C_W(E)=C_G(E)$, the result follows.
Otherwise, $[C_G(E):C_W(E)]=2$.  By Lemma \ref{abelian in sym},
all characters of
$C_W(E)$ are rational and so by Lemma \ref{index 2} and the result for
$C_G(E)$, all representations of $C_W(E)$ are defined over $\QQ$.

Next suppose that $W$ is a Weyl group of type $B_n$.  Then $W$ is
the centralizer of a fixed point free involution $z \in G :=\mathrm{Sym}_{2n}$.
We may as well assume that $z \in E$.   Then $C_G(E)=C_W(E)$ and so the result
holds.

If $W$ has type $D_n$ and $n$ is odd, then $W \times Z \cong B_n$ where
$Z$ is the center of $B_n$ of order $2$ and so the result follows by the
result for type $B$. If $W$ has type $D_n$ and $n$ is even, then
$W$ contains the central involution $z$ of $B_n$ and again we may assume
that $z \in E$.   Let $H=\mathrm{Alt}_{2n}$.  So $W= B_n \cap H$
and so $C_W(E)=C_H(E)$ and so the result follows by the result for
alternating groups.

If $W$ is exceptional, this is proved in \cite{gg}.

\end{proof}

For the other real reflection groups, we know the following:

\begin{theorem}  Let $W$ be a finite real reflection group and
$E$ an elementary abelian $2$-subgroup of $W$.
 All representations of $C_W(x)$ are defined over $\RR$.
\end{theorem}

\begin{proof}  Clearly, we may assume that $W$ is irreducible.

If $W$ is a  Weyl group, this follows
by the previous result. 

If $W$ is an irreducible real reflection
group other than a Weyl group, this follows by inspection (these are just
dihedral groups, $H_3$ and $H_4$), whence (1) follows.
\end{proof}


\section{Some Real Groups} \label{realgroupsec}

Recall that a Young subgroup of $S_m$ is the largest
subgroup of $S_m$ with a given set of orbits
(so of the form a direct product of symmetric groups
of size $m_i$ with $\sum m_i = m$).

\begin{lemma} \label{real groups} Let $A$ be a finite abelian
group of exponent $n$  and
$m$ a positive integer. Set $B=A^{(m)}$.  Let  $S_m$
act  on $B$ in the obvious way (permuting the $m$ copies
of $A$).
Let $\tau$ be the involution acting as inversion
on $B$.  Note that $\tau$ commutes with the action of
$S_m$ on $B$.   Set $G =B\rtimes(S_m\times \langle \tau \rangle)$
and $G_0=B \rtimes S_m$.
Let $L$ be the real subfield  of $K:=\QQ[\zeta]$ where
 $\zeta$ is a primitive $n$th root of unity.
Then every representation of $G$ is defined over $K$
and every representation of $G_0$ is defined over $L$.
\end{lemma}

\begin{proof}   We leave the second statement as an exercise
(it is essentially the first case  of the proof of the first statement
given below).

 Let $V$ be an
irreducible $G$-module.  Let $V_{\gamma}$ be the $\gamma$
eigenspace for $B$ for some $\gamma \in B^*$, the dual group of $B$,
and assume that this is nonzero.  Then $G$ permutes these eigenspaces
and the stabilizer $S$ of this eigenspace satisfies one
of the following possibilities:

\begin{enumerate}
\item $S=BY$, where $Y$ is a Young subgroup;
\item $\gamma^2=1$ (or equivalently, $\gamma=\gamma^*$)
and $S=B(Y \times \langle \tau \rangle)$; or
\item $S=B(\langle Y, g\tau \rangle)$
where $g \in S_m$ is an involution and $\gamma^g=\gamma^*$.
\end{enumerate}

The first case is where $\gamma$ and $\gamma^*$ are not
in the same orbit under $S_m$, the third case is
where $\gamma \ne \gamma^*$ but $\gamma$ and $\gamma^*$
are distinct but in the same orbit (and we can always choose an involution
which takes $\gamma$ to $\gamma^*$) and the second
case is where $\gamma$ is self dual.

Consider the first case.  Let $U=V_{\gamma} \oplus V_{\gamma^*}
= V_{\gamma} \oplus \tau V_{\gamma}$.  Then the stabilizer
of $U$ is $T:=B(Y \times \langle \tau \rangle)$.  By construction,
the image $D$ of $B\langle \tau \rangle$ acting on $U$ is
a dihedral group.  Thus, the image of $T$ factors through
$D \times Y \times \ZZ/2$.  Since every representation
of this group is defined over $L$, it follows that $U$ is defined
over $L$.  Since $V$ is induced from $U$, the same is
true for $V$.

The second case is similar.  Let $U=V_{\gamma}$.
Thus $S$ acts on $V_{\gamma}$
as a homomorphic image of  $Y \times \ZZ/2 \times \ZZ/2$,
and the representation is defined over $\QQ$.

Next consider the third case.  Write $\gamma=(\gamma_1, \ldots,
\gamma_m)$.
The hypothesis implies that for each $\alpha$ in the dual group $A^*$ of $A$,
the number of $i$ with $\gamma_i=\alpha$ is the same as the number
of $i$ with $\gamma_i=\alpha^*$.  It follows that $S$ has a normal
subgroup that is
the direct product of symmetric groups and that $g$ and $\tau$ are
commuting
involutions with $g\tau$
an involution which either commutes with a direct factor or interchanges
two direct factors (depending upon whether the corresponding character
is self dual or not).

Let $U = V_{\gamma}
\oplus V_{\gamma^*}$.  The stabilizer of $U$
is $T:=\langle S, \tau \rangle$.
Note that $g\tau$ commutes with the image of $B$ on $U$ and
so $T$ acts on $U$ as a homomorphic image of $D \times \langle Y,
g\tau \rangle$.  Thus, it suffices to prove that $\langle Y, g \tau \rangle$
has all representations defined over $F$.    Note that $Y$ is a direct
product of symmetric groups and $g\tau$ either commutes with a direct
factor or interchanges two factors.   So $\langle Y, g\tau \rangle$
contains a direct product $Y_1 \times \ldots Y_t$ where each
$Y_i$ is either a symmetric group a direct product of $2$ symmetric
groups and is contained in a direct product of $Z_i, 1 \le i \le t$
where $Z_i  \cong Y_i \times \ZZ/2$ or is a wreath product
$U \wr \ZZ/2$ where $U$ is a symmetric group.   Since all representations
of $Y_i$ and $Z_i$ are defined over $\QQ$, it is easy to see that
$\langle Y, \tau \rangle$ has all characters rational.
Since the representations
of $Y_i$ are defined over $\QQ$, the same is true for
$\langle Y, \tau \rangle$
by Lemma \ref{index 2}.   This completes the proof.
\end{proof}

We will need the following result to deduce the results for Weyl
groups of type D from those of type B.

\begin{lemma}\label{typeD lemma} Let $A$ be a finite cyclic group
of even order $r$ and $m$ a positive integer. Set $B=A^{(m)}$.  Then $S_m$
acts on $B$ in the obvious way (permuting the $m$ copies
of $A$).  Let $A^2$ denote the group of squares in $A$.
   Let $\tau$ be the involution acting as inversion
on $B$.  Set $G =B\rtimes(S_m\times \langle \tau \rangle)$.
Let $B_0$ be the subgroup of $B$ of index $2$ that is
normal in $G$
(i.e. $B_0=\{(a_1, \ldots, a_m) \in B | a_1\cdots a_m \in A^2\}$).
Let $G_0=B_0(S_m \times \langle \tau \rangle)$.  Let $L$
be the real subfield of $K:=\QQ[\zeta_r]$.
Then
\begin{enumerate}
\item every representation of $B_0 \rtimes S_m$ is defined over $K$;
\item every representation of $G_0$ is defined over $L$; and
\item  $\tau$ inverts all conjugacy classes in $B_0S_m$.
\end{enumerate}
\end{lemma}

\begin{proof}   View $A$ as subgroup of $\sym_{rm}$ acting
semiregularly.   Then
$B_0S_m$ is isomorphic to the centralizer of $A$ in
$\alt_{rm}$ and $G$ is the subgroup of $\alt_{rm}$
consisting of all elements that either centralize $A$ or act
as inversion on $A$.  We have already proved (3) in this context.

It is an easy exercise to show that if $r$ is $1$ modulo $m$
and relatively prime to the order of an element $g \in BS_m$,
then $g^r$ and $g$ are conjugate via $B_0S_m$ (indeed,
this follows by our result for abelian subgroups of the Weyl
group of type D).

Using Lemma \ref{index 2} and the previous result implies that
 $B_0S_m$ has all representations defined over
$K$.

In order to prove the second statement, by Lemma \ref{index 2} and
the previous result, it
suffices to prove that every representation of $G_0$ has character values
in $F$.   Since all the character values are in
$K$ by the previous paragraph,
it suffices to show that all character values are real.

Since $\tau$ inverts all conjugacy classes in $B_0S_m$,   we need
only show that if $g=bs\tau$ with $b \in B_0$
and $s \in S_m$, then $g$ and $g^{-1}$ are conjugate
in $G_0$.  We prove the formally stronger
statement allowing $b \in B$ rather and show $g$ is conjugate
to $g^{-1}$ via an element in $G_0$.

By induction on $m$, it suffices to prove this result when $s$ is
an $m$-cycle.   Note that we can replace $g$ by any $G$-conjugate.

There is no harm in assuming that $s$ is the
$m$-cycle $(1 \ 2 \  \ldots \ m)$.  By conjugating by an element
of $B$, we may assume that $b=(a_1, \ldots, a_m)$ where
$a_i \in A^2$ for $i > 1$.   If $a_1 \in A^2$, then
$g \in (A^2)^{m}(S_m \times \langle \tau \rangle$ and all representations
of this group are defined over $K$ by Lemma \ref{real groups}.
If $a_1$ is not a square
in $A$, then $g$ is not in $B_0$.  Since $g^x=g^{-1}$ for some
$x \in G$, one of $x$ and $gx$ is in $G_0$ and conjugates $g$ to
$g^{-1}$.
 \end{proof}

\section{Classical Weyl Groups }\label{typeB}

Let $G$ be a classical Weyl group.   We prove Theorem \ref{main}
and (1) of Theorem \ref{Ng} for the classical Weyl groups in this section.
By Corollary \ref{mod red for D(G)}, it suffices to work over
$\BC$.

We have already shown that
all Schur indicators for $C_g$ are $+1$ when $g^2=1$
(and all irreducible modules are defined over the prime field).

Now consider $g^2 \ne 1$.  We have already shown that
if $x \in N_g \setminus{C_g}$, then $x$ inverts all conjugacy
classes of $C_g$.
So all that remains to prove is that all Schur indicators are $+1$
for $N_g$ when $g^2 \ne 1$.

We first handle the case of $\Sym_n$.  We prove a slightly more general
result.  Namely, let $J$ be an abelian subgroup of $\Sym_n$ with centralizer
$C$.  Let $N=\langle C, \tau \rangle$ where $\tau$ acts as inversion on $C$.
If $J$ is elementary abelian, then we have shown that all representations
of $C = N$ are defined over $\QQ$.  So assume that $J$ has exponent
greater than $2$.  If $J$ has an orbit that is not regular,
then $C = C_1 \times \ldots \times C_t$ where $C_i$ is the centralizer
of $J$ acting $\Omega_i$ where $\Omega_i$ is the union of the $J$-orbits
with given point stabilizer.  Applying Lemma \ref{dp1} yields the result by
induction.  So we may assume that all orbits of $J$ are regular.
Then $N$ is the group described in Lemma \ref{real groups} and we have seen
that all representations are defined over $\RR$.

The case where $G$ is of type $B_n$  now follows (for the more general
case of abelian subgroups $J$ containing $Z(G)$) since $C$ and $N$
will be contained in $G \le \sym_{2n}$.

Finally, consider the case that $G=D_n$.  By Lemma \ref{index 2}
and the result for $B_n$, it suffices to prove that all irreducible
representations of $N_g$ are self dual.  It follows by
Theorem \ref{inversion} that $y^x$ is conjugate to $y^{-1}$ in $N_g$ for
any $y \in C_g$ and any $x \in N_g \setminus{C_g}$.

So it suffices to show:

    \begin{lemma}  Let $W$ be a Weyl group of type A, B or D.
    Let $g \in W$ and $h \in N_g \setminus{C_g}$.  Then there exists
    $x \in N_g$ inverting $h$.
    \end{lemma}

    \begin{proof}  In fact, since we have proved that all representations
    of $N_g$ have Schur indicator $+1$ for $W$ of type A and B, the result
    holds in that case.  We prove a slightly stronger result for type B
    that will give the result for type D as well.

    Let $W = B_n$.  Let $W_0$ be the subgroup of
    index 2 of type D in $W$. We prove the slightly stronger result
    that given $h  \in N_g \setminus{C_g}$,
    there exists  $x \in W_0 \cap N_g$.  There is also no harm in
    assuming that $g$ has even order $r$ (for if $g$ has odd order,
    replace $g$ by $gz$ where $z$ is the central involution of $W$).

    Since all representations of $N_g$ are defined over $\RR$, it follows
    that there is some $x \in N_g$ inverting $h$ (indeed, by replacing
    $x$ by $xh$, we can choose $x \in C_g$).  We need to show that
    we can choose $x \in C_g \cap W_0$.

    We induct on $n$.  The case $n \le 3$ is clear.
    Set $X = \langle g, h \rangle$.

    If $C_W(X)$ is not contained in
    $W_0$, the result follows as well (for choose w in $C_W(X)$ with
     $w$ not in $W_0$, then either $x \in W_0$ or $xw \in W_0$ and both
    elements inverts $h$).

In particular, if $n$ is odd, the result follows for then $Z(W)$ is
not contained in $W_0$.

We view $W \le \sym_{2n}$ centralizing the fixed point free involution
$z$.
Note  that if $X$ is contained in a direct product of
subgroups of $W$ of type B, the result holds by induction.  In particular,
this implies that we may assume that $g$ is homogeneous
-- i.e. viewing $W \le \sym_{2n}$, $g$ has all orbits of size
$r = 2s$.   The number of orbits is $m:=n/s$.  Similarly,
$\langle g, z \rangle$ has all orbits of the same size $r'$.
So either $z$ is a power of $g$ and $r'=r$ or $r'=2r$.

If $C_W(g) \le W_0$, then as $g^y=g^{-1}$ for some $y \in
W_0$ (since all representations of $W_0$ are real), $N_g \le W_0$
and $x \in N_g \le W_0$. This is the case if $r'=2r$.  So
we may assume that $r'=r$, i.e. $z$ is a power of $g$.

In this  case, $C_g \cong \ZZ/2s \wr S_m$.
It follows that $N_g = A^m \rtimes (S_m \times \langle \tau \rangle)$
where $A$ is cyclic of order $2m$ and $\tau$ inverts $A^{t}$.
By Lemma \ref{typeD lemma}, $N_g$ is a real group and so the conclusion
holds.
\end{proof}


\section{General Real Reflection Groups}\label{general}

We have now shown that if $W$ is a Weyl group of classical
type and  $g \in W$, then:
\begin{enumerate}
\item $N_g$ is totally orthogonal for every $g \in W$;
\item Every element of $C_g$ is inverted by an element
of $N_g \setminus{C_g}$.
\item $D(W)$ is totally orthogonal.
\end{enumerate}

It is obvious that all three properties hold for dihedral
groups as well.  The remaining real reflection groups
are $H_3$, $H_4$ and the Weyl groups of type $F_4$, $E_6$,
$E_7$ and $E_8$ (see \cite{H}).  Since $H_3 = A_5 \times
\ZZ / 2$, it is
clear that all three conditions hold for $H_3$ as well.

For the remaining $5$ groups, we compute using
MAGMA or GAP.   All the character
tables of these groups are in MAGMA (or GAP) and it is straightforward
to loop through the irreducible modules $C_g$ and $N_g$ and compute
their indicators.  We see that over $\BC$,
$C_g$ is totally orthogonal for $g^2=1$
and $\nu(V) \ge 0$ for all irreducible $\BC N_g$-modules $V$.

Using Proposition \ref{localind} and the character tables for
$N_g$, we see that $\nu(\hat{V}) = +1$ for irreducible $D(G)$-modules
(over $\BC$).  By the reduction lemma, this holds for any odd
characteristic as well.

We also note for $H_4$ and the Weyl groups of
type $E$ and $F$ that $N_g$ need not be totally orthogonal.
Indeed, in each case, there is an element $g$ of order $3$ such
that some irreducible $\BC N_g$-module has Schur indicator $0$
(but all indicators for $\BC N_g$ are non-negative).

This suggests that the method of proof used for the classical Weyl groups
will not extend to all Weyl groups.   The example in the next
section also shows that even being split over $\mathbb{Q}$
does not imply that $D(G)$ cannot have irreducible modules
with $0$ indicator.

We end this section by proving a result
which includes Corollary \ref{formula}.

\begin{theorem} Let $G$ be a finite group with Drinfel'd
double $D(G)$ (over $\BC$).  Let $t$ be the number of
elements $g \in G$ with $g^2=1$. Then
$$
t^2 \le \sum_{g \in G} \sum_{\chi \in \mathrm{Irr}(C_g)} \chi(1),
$$
with equality if and only if $D(G)$ is totally orthogonal.
\end{theorem}

\begin{proof}  Let $S$ be the antipode of $D(G)$.  Then it follows
from \cite{LM} that the trace of $S$ on $D(G)$ is the number
of pairs $(a,g)$ with $a^2=1=(ag)^2$.  It is easy to see that
this is $t^2$ .

It also follows from \cite{LM}
that this trace is $\sum \nu(W) \dim W$ where the
sum is over all distinct irreducible $D(G)$-modules.  Let $\Omega$
be a set of representatives of conjugacy classes of $G$.

It follows from the
description of the irreducible $D(G)$-modules, this is
$$\sum_{g \in \Omega} [G:C_g] \sum_{V \in \mathrm{Irr}(C_g)} \dim V
= \sum_{g \in G}  \sum_{V \in \mathrm{Irr}(C_g)} \dim V.$$

Thus, the inequality holds with equality precisely when
$\nu(W)=+1$ for irreducible $D(G)$-modules.
\end{proof}

See \cite{JM} for a generalization of this result.


\section{An  Example} \label{example}

We now give an example to show that one cannot extend the
main theorem to groups in which all representations are
defined over $\RR$ (or even $\QQ$).

 Let $G = \langle x,a,b | x^8=a^2=b^2=1, axa=x^{-1}, bxb=x^3 \rangle$.

Note that $G$ has order $32$ and is the holomorph of the cyclic
group of order $8$ (i.e. it is the semidirect product of the cyclic
group of order $8$ and its automorphism group).

\begin{lemma}  All representations of $G$ are defined over $\QQ$.
\end{lemma}

\begin{proof}  Note that $x^4$ is central in $G$ and
$G/\langle x ^4 \rangle \cong D_8 \times \ZZ/2$, where
$D_8$ is the dihedral group of order $8$.  All representations
of $D_8$ are easily seen to be defined over $\QQ$, hence also
those representations of $G$ with $x^4$ acting trivially.

We now construct a faithful $4$-dimensional absolutely irreducible
representation of $G$ defined over $\QQ$.  Namely, let $V= \QQ (\zeta_8)$, the field
generated by the $8$th roots of unity over $\QQ$.  This is a 4-dimensional
module over $\QQ$; $G$ acts on it by letting $x$ act by multiplication by $\zeta_8$
and $\langle a,b \rangle$ by $\mathrm{Gal}(V)$. It is easy to see that these actions give a group
isomorphic to $G$.

Now $x$ has 4 distinct eigenvalues in this representation (the primitive $8^{th}$ roots of 1)
and $G$ permutes these transitively. Any submodule must be generated by
$\langle x \rangle$-eigenspaces, and therefore is either 0 or contains all the eigenspaces. Thus
$V$ is an (absolutely) irreducible 4-dimensional module defined over $\QQ$.

Adding the sum of the squares of the degrees, we see that this gives all modules and all are
defined over $\QQ$.
\end{proof}

\begin{lemma}  Let $g =xa$.  Then $g^2=1$ and
$C_G(g) = \langle x^2ab, g \rangle
\cong \ZZ/4 \times \ZZ/2$.
\end{lemma}

\begin{proof}  It is straightforward to compute that $g$
is an involution and  commutes with $x^2ab$, an element of order
$4$.    Since $C_G(g) \cap \langle x \rangle
= \langle x^4 \rangle$, it follows that $|C_G(g)| \le 8$.
since $C_G(g) \geq \langle x^2ab, g \rangle$, we have equality.
\end{proof}

In particular, we see that for the involution $g$, $C_g$ has two
non self dual representations (this can also be seen using MAGMA).
This implies the Schur indicators
for the corresponding modules for $D(G)$ have Schur indicator
$0$ as well.

So $G$ has all indicators $+1$, but $D(G)$ has some indicators
$0$.\\

{\bf Question}:  Can there be negative indicators for $D(G)$ if all
representations of $G$ are defined over $\RR$ or $\QQ$?

See \cite{gow1, gow2, tz} for examples of totally orthogonal groups.


\end{document}